\documentclass{amsart}
\usepackage{amsmath, amssymb, amsthm, latexsym}
\input xypic
\newtheorem{theorem}{Theorem}
\newtheorem{proposition}[theorem]{Proposition}

\newtheorem{lemma}[theorem]{Lemma}

\newtheorem{definition}[theorem]{Definition}
\newtheorem{remark}[theorem]{Remark}
\newtheorem{example}[theorem]{Example}

\def\Proof{\medskip\noindent{\bf Proof: }}

\def\Z{\mathbb{Z}}

\def\A{\mathbb{A}}

\def\C{\mathbb{C}}

\def\C{\mathbb{C}}

\def\N{\mathbb{N}}

\def\Pi{\mathbb{P}^{\infty}}

\def\qed{\hfill$\square$\medskip}

\def\Zpk{\mathbb{Z}/p^{k}}
\def\Zpk1{\mathbb{Z}/p^{k-1}}

\newcommand{\rref}[1]{(\ref{#1})}

\newcommand{\cform}[3]{\begin{array}{c}
{\scriptstyle #3}\\
#1\\
{\scriptstyle #2}\end{array}}

\newcommand{\beg}[2]{\begin{equation}\label{#1}#2\end{equation}}
\def\r{\rightarrow}
\def\mc{\mathcal{C}}

\def\sl2{\widetilde{SL_{2}(\Z)}}

\title{A universal approach to vertex algebras}
\author{Ruthi Hortsch, Igor Kriz and Ale\v{s} Pultr}
\thanks{The first author was supported by the NSF REU program.
The second author was supported in part by NSA grant H98230-09-1-0045. The
third author was supported by project 1M0545 of the Ministry of Education of the
Czech Republic}

\begin{document}

\maketitle

\section{Introduction}

The notion of vertex algebra due to Borcherds \cite{bor} and 
Frenkel-Lepowsky-Meurman \cite{flm} is one of the fundamental concepts 
of modern mathematics. The definition is purely algebraic,
but its central point of the
is an intricate Jacobi identity (\cite{flm}), which makes
the definition hard to motivate from first principles. The identity
can be replaced by a locality
axiom \cite{dl}, cf. \cite{frenkel,kac}, see also \rref{edvert1} below, 
which however feels more like a part of physics than algebra. Of course,
vertex algebras are an algebraic model a part of conformal field theory,
which enters into the motivation, and also suggests generalizations
(cf. Borcherds \cite{bor1}, Soibelman \cite{soib} and
Beilinson-Drinfeld \cite{beidr}).

\medskip

In this paper, we characterize vertex algebras 
by exploring the algebraic structure
present on the {\em correlation functions} of a vertex algebra. 
This turns out be a structure of a co-operad, the dual of the 
notion of operad. Operads are known to play an important
role in algebra (cf. Ginzburg-Kapranov \cite{gk}) and topology (cf. May \cite{g}).
In algebra, operads describe vector spaces
with additional multilinear operations, which are
subject to universal multilinear identities. 
Our axiomatization of vertex algebra via the co-operad
structure on correlation functions does, in a way, explain the Jacobi identity
within the confines of algebra. We stress that our use of operads is different
from the paper of Huang and Lepowsky \cite{lh}, whose operads (or partial operads)
consist of points, while our co-operads consist of functions (in a vague sense, the
approaches could be
called ``Koszul-dual'').

\medskip

When looking at the picture in more detail, subtleties emerge, 
dealing with which comprises a large part of this paper. The extant definitions
of vertex algebra (\cite{bor, flm, frenkel, kac}) vary slightly; to get a universal
algebra treatment, we need to fix notions. We end up focusing on two notions which
differ by a degree of {\em uniformity} of the locality required. In any case,
however, if we deal with general (non-connective) vertex algebras, the category
lacks the properties of a category of universal algebras, which forces us to 
use co-operads instead of operads. The algebraical and categorical discussions
resulting are in a way standard, but as far as we know, are not in the literature,
which is why we include an appendix (Section \ref{sapp}) dealing with the prerequisites.

\vspace{3mm}

When we assume connectivity, (a vertex algebra is $k$-connective when $V_n=0$ for $n<k$),
we pass to a sub-co-operad which is finite-dimensional, and thus can be 
dualized to an actual operad, which we can, at least in a certain
sense, describe explicitly (Theorem \ref{tfiltr} below). 
Thus, $k$-connective vertex algebras for a fixed $k$
are actual universal algebras, and can be defined in terms of generators and defining
relations. In Section \ref{sgdr}, we give some examples
of using this method, although 
this is just a beginning, and more work
needs to be done. We should also mention the very important work of Roitman who in a series of papers 
\cite{roit1,roit2,roit3} explored the category of vertex 
algebra restricted by specific bounds of locality between given generators. In his category,
Roitman also has a notion of construction from generators and defining relations (in particular,
he studies free objects), but his category isn't an actual category of universal algebras,
and his results are different from ours. For examples, Roitman's
free objects typically have no connectivity, which has a physical explanation
in the Coulomb gas realization of conformal field theories (cf. \cite{gh}).
Still, in certain cases, (when a Roitman presentation enjoys connectivity),
presentations in
our sense can be deduced from Roitman's (e.g. the case of positive-definite lattices,
see Section \ref{sgdr} below). Finally, we should also mention the related work of Bokut,
Fong and Ke \cite{bokut}, who studied presentations of conformal algebras, and even developed
a Gr\"{o}bner basis style algorithm in that case.

\medskip

The present paper is organized as follows: In Section \ref{s2},
we describe the co-operad of correlation functions and its $k$-connective
version. In Section \ref{s3}, we prove our main results characterizing
vertex algebras and $k$-connective vertex algebras in this setting. 
Section \ref{sgdr} contains examples of presentations of vertex algebras
in terms of generators and defining relations. Section \ref{sapp} is the
Appendix as mentioned above.

\vspace{3mm}
\noindent
{\bf Acknowledgements:} We are thankful to Yi-Zhi Huang and James Lepowsky
for discussions.

\vspace{3mm}

\section{VA-algebras and the correlation function co-operad}
\label{s2}

In this paper, all vector spaces we consider will be over the field
of complex numbers $\C$.
We shall say that a meromorphic function $f(z)$ on $\C P^1$ (the one-dimensional
complex projective space) is {\em non-singular
at $\infty$} when the limit of $f(z)$ at $\infty$ is $0$.

\begin{definition}
\label{dlocal}
{\em A {\em local function} is a meromorphic $n$-variable function on $\C P^1$
which is non-singular
when all the variables are different and different from $\infty$. }
\end{definition}

The space of all local functions in variables $z_1$,...$z_n$ will be
denoted by $\mc(z_1,...,z_n)$.

\vspace{3mm}
\noindent
{\bf Comment:} Local functions can also be characterized as regular functions
(in the algebraic sense) on the affine variety consisting of ordered $n$-tuples
of distinct points of the affine space $\A^1$ (the ``ordered configuration space'').
Most of our arguments have a parallel in the category of algebraic varieties over
any field of characteristic $0$ (we will point out where changes are needed).

\vspace{3mm}
\begin{lemma}
\label{l1}
There exists a vector space basis $B_n$ of $\mc(z_1,...,z_n)$ such that
\beg{e1l1}{B_{0}=\{1\},
}
\beg{e2l1}{B_{n+1}=\{(z_{n+1}-z_i)^k,\; z_{n+1}^{\ell}\;|\: i=1,...,n,0>k\in\Z\; 0\leq \ell\in\Z\}\cdot B_{n}.
}
\end{lemma}

\Proof
Let $f=f(z_1,...,z_{n+1})\in\mc(z_1,...,z_{n+1}).$
Then consider $f$ as a function of $z_{n+1}$ with $z_1,...,z_n$
constant. Since every holomorphic function on $\C P^1$ is constant,
there exist unique functions $g_k\ell$, $g_{1k}$,...,$g_{nk}$,
$0>k\in\Z$, $0\leq\ell\in\Z$,
all but finitely many of which are $0$,
such that
\beg{e3l1}{f=\sum_k g_\ell z_{n+1}^{\ell} + \sum_{i,k} g_{ik}(z_{n+1}-z_i)^{k}.
}
Clearly, $g_\ell, \; g_ik\in \C(z_1,...,z_n)$, so an induction completes the proof.
\qed

\vspace{3mm}
Clearly, there is a right action of $\Sigma_n$ on $\mc(z_1,...,z_n)$
given, for a permutation $\sigma:\{1,..,n\}\r\{1,...,n\}$, by
\beg{eeq}{f\sigma(z_1,...,z_{n})=f(z_{\sigma(1)},...,z_{\sigma(n)}).}
We shall also denote $\mc(z_1,...,z_n)$ by $\mc(n)$, $n\geq 0$.

\vspace{3mm}
Now note that it follows from Lemma \ref{l1} that
each of the vector spaces $\mc(n)$ is graded by {\em the negative
of
the degree} (a function $f(z_1,...,z_n)$ is homogeneous of degree $k\in\Z$
if $f(az_1,..,az_n)=a^kf(z_1,...,z_n)$. Denote by $\mc(n)_k$ the subspace
of elements of degree $k$. (The sign reversal will be needed when
relating our concept to the concept of vertex algebras; in
writing fields of vertex algebras, the degree of a map attached as
a coefficient at $z^n$ {\em minus} $n$ is a constant.)

For $k=p+q$, we shall now describe a map
\beg{ephi}{\Phi_{n,m}:\mc(n)_k\r \mc(m+1)_p\otimes \mc(n-m)_q.}
This is done as follows: consider a function $f\in \mc(z_1,...,z_n)_k$.
By Lemma \ref{l1}, $f$ can be expressed (non-uniquely)
as a polynomial in 
\beg{emonij}{(z_i-z_j)^{-1},} 
\beg{emoni}{z_i.}
Now modify this polynomial as follows: Replace every monomial
\rref{emonij} with $i>m$, $j\leq m$ by 
\beg{emonijt}{((t-z_j)+(z_i-t))^{-1},}
and expand in increasing powers of $z_i-t$. Also 
replace every monomial \rref{emoni} with
$i>m$ by
\beg{emonit}{(t+(z_i-t)).}
Then substitute $t_s$ for $z_{s+m}-t$ with $s=1,...,n-m$
and $z_{m+1}$ for $t$. In particular, a term \rref{emonij}
will remain unchanged for $i,j\leq m$, while for $i,j>m$, the term
\rref{emonij} will turn into
$$(t_{i-m}-t_{j-m})^{-1}.$$
Picking out terms of degree $p$ in
the variables $z_1,...,z_{m+1}$, and degree $q$ in the variables
$t_1,...,t_{n-m}$,
we obtain an element
\beg{ephip}{g\in \mc(z_1,...,z_{m+1})_p\otimes 
\mc(t_1,...,t_{n-m})_q
}
(one easily sees that the number of summands of the given fixed
degree is finite). Put $\Phi_{m,n}(f):=g$.

\vspace{3mm}

\begin{lemma}
\label{econphi}
The element $g$ of \rref{ephip} does not depend on the representation
of $f$ as a polynomial in \rref{emonij}, \rref{emoni}.
\end{lemma}

\Proof
Interpret $z_1,...,z_m$ as constants different from each other and
different from $\infty$. Pick $t\in \C P^1$ different from
$z_1,...,z_m$ and from $\infty$. Then, putting
$$P=\prod_{i\neq j} (z_i-z_j)^N$$
for a sufficiently large integer $N$,
$gP$ can be characterized
as the appropriate homogeneous summand of
the Taylor expansion of $fP$ as a function in $z_{m+1},...,z_n$
in the neighborhood of 
$$(z_{m+1},...,z_n)=(t,...,t).$$
\qed

\vspace{3mm}
From now on, we will be using notation and certain facts from
category theory and the theory of operads. These facts are somewhat
technical and independent of the rest of the material, and are
treated in the Appendix. 

\vspace{3mm}
\begin{theorem}
\label{t1}
There exists a unique graded co-operad structure on the sequence
$(\mc(n))$ (we call this the {\em correlation function co-operad})
such that $\Sigma_n$-equivariance is given by \rref{eeq}
and the co-composition operation corresponding to inserting
$n-m$ variables to the $m+1$'st variable
among variables indexed $1,...,m,m+1$ is given by $\Phi_{n,m}$.
\end{theorem}

\Proof
%
The axioms we then need to verify are the version of
equivariance for insertions, commutativity of insertions into 
two different original variables, and associativity of insertions
(arising when an insertion is followed by another insertion
into the new variables).

Regarding proving these properties, first note that equivariance
is obvious, expressing simply symmetry of the construction
in the variables involved in equal capacity. Regarding
commutativity of insertions, using
\rref{emonijt}, it corresponds to the fact that
expanding
$$((s-t)+ (z_i-s)-(z_j-t))^{-1}$$
in increasing powers of $(z_i-s)$ and then increasing powers of
$(z_j-t)$ gives the same series as expanding in $(z_j-t)$ first
and in $(z_i-s)$ second. One readily verifies that the result
in both cases is
\beg{esymexp}{\sum_{m,n\geq 0} {m+n \choose n} (-1)^{n}
(z_j-s)^n (z_i-t)^m (s-t)^{-m-n-1}.}
Coassociativity is equivalent to saying that
expanding 
$$(-(z_i-s)+(z_j-s))^{-1}$$
in increasing powers of the second summand, and then 
write $(z_j-s)$ as
$$((t-s)+(z_j-t))$$
and expand in increasing powers of the second summand
gives the same result as expanding 
$$(-(z_i-t) +(z_j-t))^{-1}$$
in the second summand and then writing $-(z_i-t)$ as
$$(-(z_i-s)+(t-s))$$
and expanding in increasing powers of the second summand.
Computation shows that both computation yield
\beg{eexpc}{-\sum_{m,n\geq 0} {m+n\choose m} (z_j-t)^n(t-s)^m (z_i-s)^{-n-m-1}.}
\qed

\vspace{3mm}
\noindent
{\bf Comment:}
In our treatment of the co-operad $(\mc(n))$, we used the version of
the definition mentioned in Remark \ref{rwaste}, which 
simplifies notation considerably. However, it is not
difficult to describe the general co-operation
from Definition \ref{dop}: This is a map of
the form
\beg{egeneral+}{\begin{array}{l}\mc(z_{11},...,z_{1n_1},...,z_{k1},...,z_{kn_k})_\ell\\
\r\mc(z_1,...,z_k)_{\ell_0}\otimes\mc(t_{11},...,t_{1n_1})_{\ell_1}
\otimes ...\otimes \mc(t_{k1},...,t_{kn_k})_{\ell_k}\end{array}
}
with $\ell_0+...+\ell_k=\ell$. To get this map, start with the variables $z_{ij}$,
then introduce new variables $z_i$, $t_{ij}$, and replace $z_{ij}$
with
$$t_{ij}+z_i.$$
The terms which require expansion are 
\beg{egeneral++}{(t_{ij}+z_i-t_{i^\prime j^\prime}-z_{i^\prime})^{-1}, \; i\neq i^\prime.
}
We rewrite \rref{egeneral++} as
$$(t_{ij}+(z_i-z_{i\prime}) -t_{i^\prime j^\prime})^{-1},$$
expanding in increasing powers of $t_{ij}$ and $t_{i^\prime j^\prime}$. 
The same formulas \rref{esymexp}, \rref{eexpc} are involved in showing that this
is well defined and associative.

\vspace{3mm}
\begin{definition}
\label{dva}
{\em A VA-algebra is a graded algebra over the graded
co-operad $\mc$.}
\end{definition}

\vspace{3mm}
We shall also be interested in cutting off $VA$-algebras by
connectivity.  

\vspace{3mm}
\begin{definition}
\label{dbva}
{\em A {\em bounded VA-algebra} is a VA-algebra $X$ such that
for every $x_1,...,x_n\in X$, the image of 
$$x_1\otimes...\otimes x_n$$
under the coalgebra structure map in
$$X\otimes \mc(n)_k$$
is $0$ for $k<-N$ where $N$ depends only on $x_1,...,x_n$.}
\end{definition}

\vspace{3mm}

Now let $k\in \Z$. Consider $\mc$ as a $\Z$-sorted
co-operad in the obvious way (see Remark \ref{rconc}). 
Then in particular,
\beg{emc0}{\parbox{3.5in}{$\mc(n,n_1,...,n_m)=\mc(m)_{n_1+...+n_m-n}$ 
only depends on $n_1+...+n_m-n$}
}
(the difference $n_1+...+n_m-n$ being the degree). 

Now consider the $\coprod_{n\geq 0} \Z^{n+1}$-sorted vector
subspace $\mc^{\prime}_{k}$ of $\mc$ which in degree $(n,n_1,...,n_m)$ is
equal to $\mc(n,n_1,...,n_m)$ when $n\geq k$, and
$0$ when $n<k$. By Proposition \ref{pop}, 
Lemma \ref{lfactor} and Remark \ref{rfactor},
there is a universal largest $\Z$-sorted sub-co-operad 
$\mc_k\subset \mc$ contained in $\mc^{\prime}_{k}$.

\vspace{3mm}
\begin{definition}
\label{dkva}
{\em A $k$-connective VA-algebra is an algebra over the co-operad
$\mc_k$.}
\end{definition}

\vspace{3mm}
\begin{lemma}
\label{lfin}
For every $k\in\Z$, and every $n,n_1,...,n_m$,
the vector space $\mc_k(n,n_1,...,n_m)$ is
finite-dimensional.
\end{lemma}

\Proof
Suppose the dimension of $\mc_k(n,n_1,...,n_m)$ is
infinite. Then the degree of singularity of elements at
$z_i\r z_j$ is unbounded for some $i\neq j$ (in fact,
we may choose $i=1$, $j=2$ by symmetry).
Therefore, by Lemma \ref{l1}, we must have elements which
map by the comultiplication into non-zero elements of 
\beg{elow}{\mc(n,i,n_3,...,n_m)\otimes \mc(i,n_1,n_2)
}
with $i$ arbitrarily low. But for $i<k$, 
\beg{elow1}{\mc_k(n,i,n_3,...,n_m)\otimes \mc_k(i,n_1,n_2)=0
}
(since the second factor is $0$), which is a contradiction.
\qed

\vspace{3mm}

By Lemma \ref{lfin}, and Remark \ref{rfail}, the category of
$k$-connective VA-algebras is the category of algebras
over an operad, and in particular has free objects and
factorization (see Remark \ref{rfactor}).

\vspace{3mm}

\section{The main results}
\label{s3}

We begin with a more precise version of Lemma \ref{lfin}, which is
more technical, but is proved by exactly the same method.
Let $1\leq i_1<...<i_p\leq m$.
Let us introduce an increasing filtration on $\mc(z_1,...,z_m)$ as
follows: For $N<0$, let $F(i_1,...,i_p)_N\mc(z_1,...,z_m)=0$. For $N\geq 0$,
call $f\in \mc(z_1,...,z_n)$ {\em $N$-regular} if the following
condition holds:
Fix a set $S$ of points $z_j\in\C$, $j\neq i_1,...,i_p$. 
Assume the $z_j$'s are all different. Then 
$\prod (z_{i_q}-z_{i_s})^{k_{qs}}f$ is holomorphic on $(\C-S)^p$ for some
integers $k_{qs}\geq 0$, $\sum k_{ij}\leq N$.
 , let
\beg{efiltr}{
F(i_1,...,i_p)_N\mc(z_1,...,z_n):=
\langle f\in\mc(z_1,...,z_n)\;|\;
\text{$f$ is $N$-regular}\rangle.
}

\begin{lemma}
\label{lfiltr}
For $N\geq 0$, there exists a vector space basis $B(i_1,...,i_p)_{n}^{N}$ of 
$$F(i_1,...,i_p)_N\mc(z_1,...,z_n)$$ 
such that
\beg{e1ll1}{B(i_1,...,i_p)_{0}^{N}=\{1\},
}
\beg{e2lf}{\begin{array}{l}
B(i_1,...,i_p)_{n+1}^{N}=
\{(z_{n+1}-z_i)^k,\; z_{n+1}^{\ell}\;|\: 
\\i=1,...,n,0>k\in\Z\; 0\leq \ell\in\Z\}\cdot B(i_1,...,i_p)_{n}^{N}.
\end{array}
}
when $n+1\neq i_p$, and 
\beg{e2lf1}{\begin{array}{l}B(i_1,...,i_p)_{n+1}^{N}=
\\ \cform{\bigcup}{k+M\leq N}{}
\{(z_{n+1}-z_{i_s})^k
\;|\: s=1,...,p-1,0>k\in\Z\; \}\cdot B(i_1,...,i_{p-1})_{n}^{M}.\\
\cup \{(z_{n+1}-z_i)^k,z_{n+1}^{\ell},\; 0>k, i\neq i_1,...,i_{p-1}, 0\leq \ell\in\Z\}
\cdot B(i_1,...,i_{p-1})_{n}^{N}\end{array}
}
when $i_p=n+1$.

\end{lemma}

\vspace{3mm}

\begin{theorem}
\label{tfiltr}
If $k<0$, then $\mc_k(n, n_1,...,n_m))=0$ for all $n, n_1,...,n_m$. In  general,
$\mc_k(n, n_1,...,n_m)\subseteq \mc(z_1,...,z_m)$ is the subspace $W_k$ of all $f$
such that for any $1\leq i_1<...<i_p\leq m$, 
\beg{etflt1}{f\in F(i_1,...,i_p)_N\mc(z_1,...,z_m)
}
where
\beg{etflt2}{N=-k+n_{i_1}+...+n_{i_p}.}
\end{theorem}

\Proof
The first statement immediately follows from the second. To prove
the second statement,
our definition of the
insertion operation clearly implies that $\mc_k(n,n_1,...,n_m)$
is contained in $W_k$.
On the other hand, since the insertion operation does not
decrease the powers of the differences of the inserted variables
involved (i.e. ``does not increase singularity in the inserted variables"),
the $W_k$'s form a $\Z$-sorted sub-co-operad of $\mc$.
\qed

\vspace{3mm}

Because of variations in the definition, let us define our
notion of vertex algebra.

\begin{definition}

\label{dverta}

{\em A {\em vertex algebra} is a $Z$-graded vector space $V$
together with, for each $a\in V_k$, a series 
\beg{edvert1}{a(z)=\sum_n a(n)z^{-n-1}
}
where 
$$a(n):V_m\r V_{m-n+k-1},$$
an operator $L_{-1}:V_n\r V_{n+1}$,
and an element $1\in V_0$
such that the following properties hold:
\beg{edloc}{\parbox{3.5in}{For every $a,b,c\in V$
and every $m\in\Z$, there exists an $N\in \N$ such
that the
$V_m[[z,t]]$-summand $s_m(z,t)$ of $(a(z)b(t)-b(t)a(z))c$
satisfies $(z-t)^Ns_m(z,t)=0$.}
}
\beg{edvertd}{[L_{-1},a(z)]=\partial_{z}a(z)}
(the symbol $\partial_z$ denotes (partial) derivative by $z$),
\beg{edvertl1}{L_{-1}1=0,}
\beg{edvertu1}{1(z)=1,}
\beg{edvertu2}{a(n)1=0 \;\text{for $n\geq 0$.}}}
We will call a vertex algebra {\em uniform} when the axiom
\rref{edloc} is replaced by the following obviously stronger
statement:
\beg{edvertl}{\parbox{3.5in}{For $a,b\in V$, $a(z)$ and $b(t)$
are local in the sense that $(z-t)^N(a(z)b(t)-b(t)a(z))=0$ for
some integer $N\geq 0$ (possibly dependent on $a$, $b$).}
}
\end{definition}

\vspace{3mm}
\noindent
{\bf Comment:}
The axioms \rref{edloc}, \rref{edvertl} are known as locality axioms.
Traditionally (cf. \cite{kac}), the stronger axiom \rref{edvertl}
are used.
In this paper, we consider {\em $\Z$-graded} vertex
algebras. In this context, the axiom \rref{edvertl} has
the non-uniform version \rref{edloc}. On the other hand,
in \cite{kac}, one considers vertex algebras without a $\Z$-grading,
in which case the non-uniform axiom does not appear to make sense.
One has the following result:

\vspace{3mm}

\begin{lemma}
\cite{roit1}
\label{lbound}
A uniform vertex algebra satisfies the following property:
\beg{edvertb}{\parbox{3.5in}{For any $a,b\in V$ there exists
$K\in\Z$ such that $a(i)b=0$ for $i>K$.}
}
\end{lemma}

\Proof
By Kac \cite{kac}, a uniform vertex algebra satisfies
the identity
\beg{ejjj}{a(m)b(k)-b(k)a(m)=\sum_{j\geq 0} {m\choose j} (a(j)b)(m+k-j).
}
When \rref{edvertl} holds with a uniform $N$, then we claim that
\rref{edvertb} must hold with $K=N$. Indeed, if $a(i)b\neq 0$
for $i>N$, then for $k+m-i=-1$, 
$$a(z)b(t)-b(t)a(z)$$
applied to $1$ has the summand
$$\sum_{m\in \Z} {m\choose i} a(i)bz^mt^{-1-m+i},
$$
which is not annihilated by $(z-t)^N$.
\qed

\vspace{3mm}

\begin{theorem}
\label{tvert}
There is a canonical equivalence between the category of vertex algebras
(resp. uniform vertex algebras, resp.
vertex algebras $V$ with $V_n=0$ for $n<k$) and VA-algebras
(resp. bounded VA-algebras
resp. $k$-connective VA-algebras). Moreover, the equivalence
commutes with the forgetful functors to $\Z$-graded vector spaces.
\end{theorem}

\vspace{3mm}
\noindent
{\bf Comment:} By ``canonical'' we mean that there is an obvious preferred
choice for the equivalence (rather than just that an equivalence exists).

\vspace{3mm}

\Proof
First, assume $V$ is a vertex algebra according to the Definition
\ref{dverta}. Then generalizing the usual arguments (cf. \cite{kac,frenkel}),
one can show
that $V$ possesses 
correlation functions. In one way, an $n+2$-point correlation function
can be described, for any integers $k_0,k_1,...,k_n, k_\infty$
as a linear combination 
\beg{etvert1}{\phi(a_0,a_1,....,a_n)(0,z_1,...,z_n)=\sum b_j f_j(0,z_1,...,z_n)}
with coefficients
$b_j\in V_{k_\infty}$ of local functions (in the sense of definition 
\ref{dlocal}) in variables $z_0$, $z_1$,...,$z_n$, dependent $n$-linearly
on homogeneous elements $a_i\in V_{k_i}$, $i=0,...,n$ ,
such that
$\phi$
has expansion 
\beg{etvert2}{a_n(z_n)...a_1(z_1)(a_0)
}
convergent in the range 
$$|z_n|>...>|z_1|>0.$$
Additionally, one can choose
\beg{etvert3}{f(z_0,z_1,...,z_n)=\exp(-z_0L_{-1})f(0,z_1-z_0,...,z_n-z_0).
} 
The structure map 
\beg{evertstr}{V_{k_0}\otimes V_{k_1}\otimes...\otimes V_{k_n}
\r \mc(z_0,z_1,...,z_n)_{k_\infty-k_0-...-k_n}\otimes V_{k_\infty}
}
is then defined by
\beg{etvert4}{a_0\otimes a_1\otimes...\otimes a_n
\mapsto \phi(a_0,a_1,....,a_n)(z_0,z_1,...,z_n):=\sum b_j f_j(z_0,...,z_n).
}
Next, we show that this defines a co-operad algebra structure. First,
equivariance
is a well known property (see e.g. Kac \cite{kac}).
To prove compatibility with the insertion operator, one can use
the well known fact that 
$$\phi(a_0,a_1,....,a_n)(z_0,z_1,...,z_n)$$
also has an expansion
$$(a_n(z_n-t)...a_{m+1}(z_{m+1}-t)a_m)(t)a_{m-1}(z_{m-1})...a_{1}(z_1)a_{0}(z_0)1$$
convergent for
$$\begin{array}{l}1>|t|>|z_{m-1}|>...>|z_1|>|z_0|,\\
|t|-|z_{m-1}|>|z_n-t|>...>|z_{m+1}-t|>0.
\end{array}$$
In the uniform case, 
to prove boundedness, use Lemma \ref{l1}, Lemma \ref{lbound} and axiom \rref{edvertb} repeatedly.

On the other hand,
assuming $V$ has the structure of a graded $\mc$-algebra,
and letting, under \rref{evertstr}, for $n=1$, as above,
\beg{etvert5}{a_0\otimes a_1\mapsto \phi(a_0,a_1)(z_0,z_1),
}
we let the degree $k_\infty$ summand of
$$a_1(z)(a_0)$$ 
be the Taylor expansion at $0$ of
$$\phi(a_0,a_1)(0,z).$$
Let also, as usual, $L_{-1}a$ be defined as the constant term of
$$(a(z)1)^\prime.$$
To prove locality, first note that by the fact that a meromorphic function
on $\C P^1$ is the sum of the singular parts of its expansions at all
points (this,  in turn, is due to the fact that a non-singular function
on all of $\C P^1$ is $0$). Now the co-composition axiom and
the identity
\beg{etvert6}{z^m=\sum_{i\geq 0} {m \choose i} (z-t)^i t^{m-i}\;\text{for $|z-t|<|t|$}
}
imply that 
\beg{etvertzm}{z^m\phi(?,a,b)(0,z,t)
}
is anywhere on $\C P^1$ equal to
\beg{etvertzm1}{\begin{array}{l}\displaystyle
\sum_{p\geq m} b(t)a(p)z^{m-p-1} +\sum_{p<m}a(p)b(t)z^{m-p-1}\\[3.5ex]\displaystyle
+\sum_{0\leq i\leq j} {m \choose i} (z-t)^{i-j-1}(a(j)b)(t)t^{m-i}.
\end{array}
}
Subtracting $z$ times
\rref{etvertzm} from the analogous expression obtained by replacing
$m$ by $m+1$, we obtain
\beg{etvertlhs}{a(m)b(t)-b(t)a(m)
}
on the left hand side and
\beg{etvertrhs}{\begin{array}{l}
\displaystyle
\sum_{0\leq i\leq j} {m+1\choose i} (z-t)^{i-j-1}(a(j)b)(t)t^{m+1-i}
\\[3.5ex]
\displaystyle
-((z-t)+t)\sum_{0\leq i\leq j} (z-t)^{i-j-1} (a(j)b)(t) t^{m-i}.
\end{array}
}
Using
$${m+1 \choose i} = {m \choose i} + {m \choose i-1} \;\text{for $i>0$,}$$
after chain cancellation, the only surviving term of \rref{etvertrhs}
is the term corresponding
to ${m\choose i-1}$ and $i=j+1$ in the first summand , i.e.
\beg{etvertr}{\sum_{0\leq i} {m\choose j} (a(j)b)(t)t^{m-j}
}
Equating \rref{etvertlhs} and \rref{etvertr} and taking the
coefficient at $t^{-k-1}$ gives
\rref{ejjj},
which is a part of the Jacobi identity. Even more importantly for us, however,
the graded co-operad algebra structure implies that when applied to
an element $c$, the summand of given degree $m$ on the right hand side
of \rref{ejjj} has $j$ bounded above, say, by $N$.
Grouping terms of \rref{ejjj} 
with $m+k$ constant,
we then see that the $V$-degree $m$ summand of
\beg{eliab}{[a(z),b(t)]c}
is equal to a sum of terms which are $V$-degree $m$ summands
of some constant times
a power of $z$ times $(a(i)b)(t)c$ times
the $i$'th derivative of the delta function $\delta(z/t)$ (say, by $z$),
for $0\leq i\leq N$. Thus, \rref{eliab} is annihilated when multiplied
by $(z-t)^N$, as required. 

One also sees clearly that in the bounded case, 
the terms of \rref{ejjj} with $j>N$ for a constant $N$ dependent on $a,b$
vanish outright, which implies uniformity.

To prove the axiom \rref{edvertd}, consider the following diagram from the
definition of algebra over a co-operad:
\beg{eddiag}{\diagram
X\otimes X\rto^\theta\dto^\theta &X\otimes \mc(2)\dto^{\theta\otimes 1}\\
X\otimes \mc(2)\rto_{1\otimes\gamma} &X\otimes\mc(1)\otimes\mc(2).
\enddiagram
}
Our plan is to consider the image of $a\otimes b$ under \rref{eddiag};
we will write the variable corresponding to $a$ (resp. $b$) under
$\theta$ as $z$ (resp. $t$). We then perform insertion to the variable
$t$, calling the new variable $u$, and we shall set $t=0$, and
take the linear term with respect to $u$. 
Taking the path through the upper right corner
of \rref{eddiag}, we obtain
\beg{edderl}{a(z)L_{-1}b,
}
while taking the path through the lower left corner, we obtain
two terms corresponding to $u$-linear terms of the expansion
$$z^n=\sum_{i\geq 0}{n\choose i}(z-u)^{n-i}u^{i},$$
namely coming from $-n(z-u)^{n-1}$ and $(z-u)^nu$.
In this order, this corresponds to the first and second term
of the sum
\beg{edderr}{-a^\prime(z)b+L_{-1}a(z)b.
}
The equality between \rref{edderl} and \rref{edderr} is
\rref{edvertd}. All the other axioms are obvious, which completes
the proof of the statement of the Theorem for general
and bounded VA-algebras.

The statement for the $k$-connective follows easily: 
clearly, if we have a $k$-connective VA algebra $X$, then
by applying 
$$\theta:X\r X\otimes \mc(1)$$
and substituting $1$ for the variable $z$ (i.e. applying
the augmentation $\mc(1)\r \mc(0)$), we get the unit, but
when starting with element of degree $<k$, the corresponding
part of the $\Z$-sorted operad $\mc_k$ is $0$, so the identity
is $0$ on those elements, proving that we have a $k$-conective
vertex algebra. 

On the other hand, starting with a $k$-connected vertex algebra,
using the $\mc$-co-alegbra structure which we already
proved, co-associativity and Theorem \ref{tfiltr}, we see that
\rref{evertstr} specializes to a map
$$V_{k_0}\otimes V_{k_1}\otimes...\otimes V_{k_n}
\r \mc_k(k_\infty,k_0,...,k_n)\otimes V_{k_\infty},$$
as required.
\qed

\vspace{3mm}

\noindent
{\bf Remark:} If we want to work in the algebraic
category (over a field of characteristic $0$), the
convergence arguments in the beginning of the proof
need to be replaced with a discussion of ``order
of expansion'', as treated, say, in \cite{flm}.

\vspace{3mm}

\section{Constructions of vertex algebras from generators and
defining relations}

\label{sgdr}

By Theorem \ref{tvert} and the remarks following Lemma \ref{lfin},
the category of $k$-connective vertex algebras is equivalent to 
a category of algebras over an operad, and the equivalence
commutes with the forgetful map, so we can speak of
free $k$-connective vertex algebras on a graded vector
space and by the remarks preceding
Lemma \ref{lfactor}, also of $k$-connective vertex algebras defined
by means of generators and defining relations. Note that
although we haven't identified an exact vector space basis
of a free $k$-connective vertex algebra $C_kX$ on a graded
vector space $X$, by \rref{eopmm},
and Theorem \rref{tfiltr}, we have a reasonably explicit presentation
of $C_kX$ in terms of generators and defining relations in
the category of graded vector spaces. The main purpose of
this section is to give the presentations of certain well known
vertex alegbras by generators and defining relations.

\vspace{3mm}

\begin{remark}
\label{rnot}
{\em Clearly, the operadic notation is awkward for the purposes of
practical applications. From Lemma \ref{l1}, (or alternately 
simply Definition \ref{dverta}) we see that
every element in a ($k$-connective) vertex algebra generated by 
certain elements can be written by finite words using the binary
operation
\beg{ebin}{a(n)b, \;n\in \Z}
and unary operation
\beg{eun}{L_{-1}a}
in the indeterminates $a$, $b$.
It is appealing from this point of view to rewrite \rref{ebin} as
\beg{ebinlie}{[a,b]_n
}
and \rref{eun} as
\beg{eunlie}{a^\prime.
}
The following Lemma then suggests that vertex algebras can
be understood as a sort of  ``deformation'' of the notion of Lie
algebras. Indeed, the Lemma for example clearly exhibits
the known fact that for a vertex algebra $V$,
$V/L_{-1}V$ is a Lie algebra with respect to
the operation $[?,?]_0$. Despite of the appeal of the notation of Lemma \ref{lnot}, 
however, the standard notation
\rref{ebin}, \rref{eun}, and especially
the field notation, is often preferable in obtaining intuition from mathematical
physics. }
\end{remark}

\vspace{3mm}

\begin{lemma}
\label{lnot}
In a vertex algebra, we have
\beg{emder}{[a^\prime,b]_n=n[a,b]_{n-1},
}
\beg{emder2}{[a,b]_{n}^{\prime}=[a^\prime, b]_n+[a,b^\prime]_n,
}
\beg{eanti}{[b,a]_m=(-1)^{m+1}\sum_{j\ge 0}\frac{1}{j!} [a,b]^{(j)}_{j+m},}
\beg{enjacobi}{
[a,[b,c]_n]_m - [b,[a,c]_m]_n = \sum_{i \ge 0} {m \choose i} [[a,b]_i,c]_{m+n-i}.
}
\end{lemma}

\Proof
\rref{enjacobi} is an immediate rewrite of \rref{ejjj}. \rref{emder}, \rref{emder2}
are usual properties of the shift in vertex algebras, proved for example in
\cite{kac}. \rref{eanti} is a consequence of those properties and
\rref{enjacobi}.
\qed

\vspace{3mm}
In figuring out the structure of a vertex algebra in terms of generators
and defining relations, it is usually key to understand first the relations
among the operations \rref{ebin} for $n\geq 0$. This corresponds, in
physical language, to figuring out the operator product expansion
(OPE) of the generators. In many cases, the Existence theorem (see \cite{fkrw},
Kac \cite{kac}, Theorem 4.5) then determines explicitly the structure
of the vertex algebra in question.

More precisely, we have:

\vspace{3mm}
\begin{proposition}
\label{egl}
Consider an $m$-connected vertex algebra given
by an ordered set of generators $B$ and
relations equating every \rref{ebin} for $a,b\in B$, $n\geq 0$, $a\leq b$
with a linear combination of derivatives of elements of $B$, then 
by Proposition \ref{enjacobi},
\beg{elievabasis}{b_{1}(n_1)....b_k(n_k)(1)
}
with $b_i\in B$, $0>n_1\geq...\geq n_k$, $n_i=n_{i+1}\Rightarrow b_i\geq b_{i+1}$
generates $V$ as a vector space. 
\end{proposition}
\qed

\begin{example}
\label{elieva}
\cite{frenkel}
{\em Let $L$ be any Lie algebra with a $2$-cocycle represented as
an invariant symmetric bilinear form $\langle?,?\rangle$. Then the Lie
vertex algebra is
the $0$-connective vertex algebra $VL$ with
generators $L$ in dimension $1$, and relations 
\beg{elieva1}{[a,b]_0=[a,b],}
\beg{elieva2}{[a,b]_1=\langle a,b\rangle\cdot 1.}
Then for dimensional reasons, no further operations \rref{ebin} with
$n\geq 0$ are possible, and hence by \rref{enjacobi}, if we denote
by $B$ an ordered basis of $L$, then \rref{elievabasis}
generates $VL$, and by the existence theorem (Theorem 4.5 of \cite{kac},) \rref{elievabasis} is
in fact a vector space basis of $VL$.}
\end{example}

\vspace{3mm}

\begin{example}
\label{evir}
{\em The Virasoro vertex algebra in the category of
$0$-connective vertex algebras has one generator $L$ of dimension $2$ and relations
\beg{evirr}{[L,L]_0=L^\prime,\; [L,L]_1=2L,\; [L,L]_2=0,\; [L,L]_3=(c/2)\cdot 1
}
where $c$ is the central charge. 
Again, this closes the OPE algebra, and therefore
again a vector space set of generators of the Virasoro vertex algebra
is \rref{elievabasis} where $B=\{L\}$. Again, this is a basis 
(cf. Frenkel \cite{frenkel}) by the Existence theorem.}
\end{example}

\vspace{3mm}
Sometimes, the OPE is not sufficient for determining the structure of 
the vertex algebra completely when there are relations between the
higher derivatives of the generators.

\vspace{3mm}

\begin{example}
\label{flm}
{\em Let $L$ be an even positive-definite lattice with $\Z/2$-valued bilinear form
$B$ such that for all $x\in L$,
\beg{elatf}{b(x,x)=\langle x,x\rangle/2.
}
Then the lattice vertex algebra $V_L$ is
\beg{estrvl}{V{L_\C}\otimes\C[L]
}
where $L_\C=L\otimes\C$ is the trivial Lie algebra with symmetric
bilinear form
coming from $L$ and $\C[L]$
is the group algebra (we denote the canonical basis elements of $\C[L]$
by $(\lambda)$, $\lambda\in L$). The $V{L_\C}$ factor is a vertex subalgebra whose
structure is specified in Example \ref{elieva}. 
More generally, we interpret a product $x\otimes (\lambda)$ in \rref{estrvl}
as $x(-1)\lambda$. Then in this notation, we have 
\beg{estrvl1}{a(0)(\lambda)=\langle a,\lambda\rangle(\lambda),\;
\; a(n)(\lambda)=0, \;n>0, \;a\in L_\C,
}
and for $x\in V{L_\C}$, 
\beg{estrvl2}{(\lambda)(z)x(-1)(\mu)=(-1)^{b(\lambda,\mu)}
(:\exp (\partial_{z}^{-1}\lambda(z)):x)(-1)(\lambda+\mu)
}
where the antiderivative is interpreted in the usual way, cf. \cite{flm}.

We can therefore take the generators of $V_L$ to be 
\beg{estrlg}{\lambda_\C,
(\lambda),} 
where $\lambda\in L$. The formulas
\rref{estrvl1}, \rref{estrvl2} certainly determine the OPE of the
generators \rref{estrlg}. More concretely, recall the convention
that we write $a_+(z)$ resp. $a_-(z)$ for the sum of terms with
non-negative resp. negative power of $z$ in $a(z)$,
then the corresponding {\em normal order} expression is
$$:a(z)b(t): = a_+(z)b(t)+b(t)a_-(z).$$
Then we have
\beg{elope1}{(\lambda)(z)(\mu)(t)=:(\lambda(z)(\mu)(t):
+(z-t)^{\langle \lambda,\mu\rangle}(\lambda+\mu)(t),
}
\beg{elope2}{\mu_\C(z)(\lambda)(t) =
:\mu_\C(z)(\lambda)(t): +\langle\lambda,\mu\rangle(z-t)^{-1}(\lambda)(t),
}
\beg{elope3}{\mu_\C(z)\lambda_\C(t)=:\mu_\C(z)\lambda_\C(t): +
\langle \lambda,\mu
\rangle (z-t)^{-2}.
}
When considering operations of the form \rref{ebinlie} with $n\geq 0$,
then the normal order expression vanishes. Therefore, clearly the
relations we get from \rref{estrlg} satisfy the assumptions of Propositiion \ref{egl}
if we choose the order in such a way that
$$\lambda_\C<(\mu)$$
for any $\lambda,\mu\in L$.

Even then, however, we see that there are more expressions \rref{elievabasis}
than a basis of \rref{estrvl}, since \rref{estrvl} allows only one label generator,
and no derivatives.  We can remedy the situation by including the relation
$$(0)=1$$
and a relation following from \rref{estrvl2} which expresses
any $\lambda(n)\mu$, $n<0$ as a product of terms of the form $\mu_\C(m)$,
$m<0$ from the left with $(\lambda+\mu)$.}
\end{example}

\vspace{3mm}

Nevertheless, in this case, there is an extremely clever clever solution
given by the following 

\begin{proposition}
\label{proit}
\cite{roit2}
Let $\Pi$ be a basis of the lattice $L$, and assume that $b(x,x)=1$
for $x\in \Pi$. Then the in the category of $0$-connective vertex algebras,
$V_L$ can be presented in generators 
$$(\lambda),\; \lambda\in \pm \Pi,$$
and relations
$$(\lambda)(n)(\mu)=0\;\text{when $n<-\langle\lambda,\mu\rangle$},
$$
$$(\lambda)(||\lambda||^2)(-\lambda)=1.$$
\end{proposition}

\Proof
Roitman \cite{roit2}, Section 10 proves that the canonical morphism 
\beg{eroitc}{\psi:V\r V_L
}
from the vertex algebra $V$
given by 
the given generators and defining relations in his category to $V_L$ is an
isomorphism. Therefore, $V$ is in particular $0$-connective. 
However, since the $0$-connectivity restriction can be viewed as
additional relations in Roitman's category, \rref{eroitc} factors
through
a morphism
\beg{eroitc2}{\tau:V\r V_{0}
}
where $V_{0}$ is the $0$-connected vertex algebra with the given
generators and defining relations. Since, in the present case, $\psi$ is iso, so is $\tau$.
\qed

\vspace{3mm}

One canonical way of always selecting relations is as follows: let $V$ be a $k$-connective
vertex algebra, $k\leq 0$. Consider then the relations
\beg{emaxpos}{\parbox{3.5in}{$a=0$ for $a\in V$ whenever 
$|b_1(n_1)...b_k(n_k)a|\leq 0$ implies $b_1(n_1)...b_k(n_k)a=0$ for $b_1,...,b_k\in V$,
$n_1,...,n_k\in \Z$.}
}

\vspace{3mm}
\begin{lemma}
\label{lmaxpos}
As a vector space, the quotient of $V$ by the relations \rref{emaxpos} is
isomorphic to $V/I$ where
\beg{elmaxpos}{\parbox{3.5in}{$I=(L_{-1})^m(b_1(n_1)...b_k(n_k)a)$ where $b_i\in V$, $n_i\in\Z$
and $a$ is as in \rref{emaxpos}, $m\geq 0$.}
}
In particular, 
\beg{elnonz}{V/I\neq 0.}
\end{lemma}

\Proof
By Lemma \ref{lfactor}, the map from $V$ to the quotient by the relations
\rref{emaxpos} is onto. Further, by its definition, $I$ is clearly contained
in the kernel. Thus, by universality it suffices to show that $V/I$ is indeed a 
vertex algebra, or that $I$ is an ideal in the category of $k$-connective
vertex algebras considered as universal algebras. This means that $I$ must
be stable under all operations where precisely one entry is in $I$ and the
other entries are in $V$. Clearly, $I$ is stable under the derivative. 
Additionally, clearly $a\in V$ and $u\in I$ together imply $a(n)u\in I$, 
by \rref{emder} and \rref{emder2}. It then follows that also $u(n)a\in I$
by \rref{eanti}. The assertion \rref{elnonz} now clearly follows
from the condition on $a$ in \rref{emaxpos}.
\qed

\vspace{3mm}
We shall refer to the ideal $I$ from Lemma \ref{lmaxpos} as the {\em maximal
positive ideal}. The quotients of Lie vertex algebras (resp. Virasoro vertex
algebras) by the maximal positive ideal $I$, when $I\neq 0$, are referred to
as WZW models (resp. minimal models). (There are variants of all these
notions involving supersymmetry, but we won't discuss them here.)
One can also prove that the lattice
vertex algebra can be characterize in this way:

\vspace{3mm}
\begin{proposition}
\label{platt}
Let $V$ be the vertex algebra with generators \rref{estrlg}
and relations for $a(n)b$ where $n\leq 0$, and $a$, $b$ are generators,
given by the corresponding summands of \rref{estrvl2}, and let $I$
be the maximal positive ideal in $V$. Then $V/I$ is the lattice algebra
corresponding to $L$.
\end{proposition}

\Proof
The key point is to show that the maximal positive ideal of $V_L$ is $0$;
then there exists a map of $0$-connective vertex algebras
$$V/I\r V_L$$
which in turn must be an iso, because its kernel would be contained in
the maximal positive ideal of $V/I$, which must be $0$.

To show that the maximal positive ideal of $V_L$ is $0$, we invoke
the well known fact \cite{flm} that $V_L$ is semisimple and its
irreducible modules correspond to elements of $L^{\prime}/L$ (where
$L^\prime$ is the dual lattice). For any semisimple vertex algebra,
the maximal positive ideal is $0$ (for example by Zhu \cite{zhu}).
\qed

\vspace{3mm}
\begin{example}
\label{exmoon}
{\em The Moonshine module $V^\natural$. 
Let $B$ be Griess' commutative non-associative algebra
\cite{griess}. Then $B$ is the weight $2$ summand of the Moonshine module
\cite{flm}, and the operation in $B$ is $[?,?]_1$.
The operation $[?,?]_2$ is $0$ and the operation $[?,?]_3$ is the pairing
invariant under the action of the Monster. We can take the free $0$-connected
vertex algebra on the generator space $B$ with these relations. Let $I$
be the maximal positive ideal in $V$. Then one
has
\beg{emoon}{V/I\cong V^\natural.
}
The reason is similar as in the lattice case: $V^\natural$ is semisimple and
has only one irreducible module, so its maximal positive ideal is $0$. 
This induces a map from the left hand side to the right hand side of 
\rref{emoon}. The map is onto because $V^\natural$ is generated
by its elements of weight $2$ (see \cite{flm}).
Any kernel would again be contained in the maximal positive
ideal of $V/I$, which is $0$.

Determining completely the OPE of the Moonshine module 
in the sense of Proposition \ref{egl} is an interesting problem.
This is because we have the operation 
\beg{emop0}{[?,?]_0.
}
In fact, as we remarked above, for any vertex algebra $V$, $V/L_{-1}V$
forms a Lie algebra under the operation \rref{emop0}. In the case of $V^\natural$,
it is an interesting problem to determine the structure of this Lie algebra $L$:
the weight $3$ part of $V^\natural$ decomposes as a representation
of the Monster into dimensions
$$21493760=21296876+196883+1.$$
The second two terms are the image of $L_{-1}$. Thus, we see that $L_3$
is the irreducible representation of dimension $21493760$ and
$L_2$ is the irreducible representation of dimension $196883$.
In particular, the Lie bracket exhibits the $21493760$-dimensional representation
as a quotient of the second exterior power of the $196883$-dimensional
representation.

If we knew the complete structure of $L$, then we would know the OPE
in the sense of  Proposition \ref{egl}; the operations $[?,?]_1$ and \rref{emop0}
satisfy the relations
$$[a,[b,c]_1]_0-[b,[a,c]_0]_1=[[a,b]_0,c]_1,$$
$$[a,[b,c]_0]_1-[b,[a,c]_1]_0=[[a,b]_0,c]_1 +[[a,b]_1,c]_0,$$
which are special cases of \rref{enjacobi}.

An easier approach to finding a complete representation of the
$0$-connective vertex algebra $V^\natural$ by generators and defining
relations seems to be from the construction of Frenkel, Lepowsky and
Meurman
\cite{flm}. They constructed $V^\natural$ as the sum of the vertex algebra
$V_{L}^{\theta}$ where $L$ is the Leech lattice and $\theta$ is the involution
$\lambda\mapsto -\lambda$ and $V_{L,T}^{\theta}$ where $V_{L,T}$ is the
twisted module of $V_L$ with respect to the involution $\theta$.
Let us focus on the untwisted sector. We can choose generators
\beg{emoong1}{\mu_\C\nu_\C,
}
\beg{emoong2}{(\lambda) +(-\lambda),
}
\beg{emoong3}{\mu_\C(\lambda)-\mu_\C(-\lambda)
}
where $\lambda,\mu,\nu\in L$. These generators do not
close under the OPE operations of Proposition \ref{egl}.
However, we can express a general element of $V_{L}^{\theta}$
as a product 
\beg{emoonprod}{b_1(n_1)...b_k(n_k)a
}
where $n_i<0$, $b_i$ are of the form \rref{emoong1} and
$a$ is of the form \rref{emoong2} or \rref{emoong3}.
Using \rref{elope3}, we can see that the elements \rref{emoong1}
close under the OPE, and hence the elements $b_i$ in \rref{emoonprod}
can be reordered using their OPE relations. Using \rref{elope2} and
\rref{elope3}, one can also write relations equating the appropriate
elements
\beg{emoong4}{b(n)a, \; n<0} 
where $b$ is of the form \rref{emoong1} and $a$ is
of the form \rref{emoong2} (corresponding to reordering the individual
$\mu_\C$ factors in \rref{elope3}. Finally, using
\rref{elope1} - \rref{elope3} again,
we can describe relations for applying \rref{ebin} with $a$ of the form
\rref{emoong2} or \rref{emoong3} where $||\lambda||=2$ to
\rref{emoong1}-\rref{emoong3} and yielding results of the
form \rref{emoonprod}. We omit the details. Putting all these 
relations together, we obtain a presentation of $V_{L}^{\theta}$
in terms of generators and defining relations. Adding relations
describing products \rref{ebin} where $a,b$ are weight $2$ elements
of $V_{T}^{\theta}$, we could write down an explicit presentation of $V^\natural$
in terms of generators and defining relations.}
\end{example}

\vspace{3mm}
\section{Appendix: Some facts on categories and operads}
\label{sapp}

Universal algebras with underlying structure of a vector space
where the additional operations are distributive can be often
axiomatized as algebras over an operad. Co-operads 
become useful when certain finiteness conditions fail.
Because of variations in the definitions, we recall here the definition
of an operad, co-operad, and universal algebra over them. 
We will work in the underlying category $Vect$ of $\C$-vector spaces and
homomorphisms.
We will also consider the category $\Z-Vect$ of $\Z$-graded $\C$-vector
spaces and homomorphisms preserving degree.

\vspace{3mm}
\begin{definition}
\label{dop}
{\em An {\em operad $\mc$} in $Vect$ is a sequence of objects $\mc(n)$, $n=0,1,2,...$
with right $\Sigma_n$-action,
and morphisms 
\beg{dopg1}{\gamma:\mc(n)\otimes \mc(m_1)\otimes...\otimes \mc(m_n)\r
\mc(m_1+...+m_n)}
such that we have equivariance, stating that for permutations $\tau\in\Sigma_n$,
$\sigma_i\in \Sigma(m_i)$,
\beg{edop1}{\diagram
\mc(n)\otimes\mc(m_1)\otimes...\otimes\mc(m_n)\dto_{Id\otimes\sigma_1
\otimes...\otimes \sigma_n}\ddrto^{\tau\wr \sigma_1,...,\sigma_n}&\\
\mc(n)\otimes\mc(m_1)\otimes...\otimes\mc(m_n)\dto_{\tau\otimes T}&\\
\mc(n)\otimes\mc(m_{\tau(1)})\otimes...\otimes\mc(m_{\tau(n)})\rto_(.6){\gamma}
& \mc(m_1+...+m_n).\enddiagram
}
($\wr$ denotes the wreath product of permutations and $T$ is the switch)
and associativity, which states that
\beg{edop2}{\diagram
\mc(n)\otimes M\otimes L_1\otimes...\otimes L_n\ddto_{Id
\otimes\gamma\otimes...\otimes\gamma}\drto^{\gamma\otimes Id\otimes...\otimes Id}&\\
&\mc(m_1+...+m_n)\otimes 
L_1\otimes...\otimes L_n
\dto_{\gamma}\\
\mc(n)\otimes\mc(\ell_1)\otimes...\otimes\mc(\ell_n)\rto_{\gamma}&\mc(\ell_1+...+\ell_n)
\enddiagram
}
where $\ell_i=\ell_{i1}+...+\ell_{ik_i}$,
$$\mathcal{M}=\mc(m_1)\otimes\mc(m_2)\otimes\mc(m_n),$$
$$L_i=\mc(\ell_{i1})...\otimes\mc(\ell_{ik_i}).$$
Since \rref{edop1}, \rref{edop2} are expressed in terms of diagrams, 
we may therefore define a {\em co-operad} as a sequence of objects $\mc(n)$,
$n=0,1,2,...$, a left $\Sigma_n$ action on $\mc(n)$ and morphisms
\beg{dopg2}{\gamma:\mc(m_1+...+m_n)\r \mc(n)\otimes \mc(m_1)\otimes ...\mc(m_n)}
such that diagrams dual to \rref{edop1} and \rref{edop2} commute. 
Definitions of operads and co-operads over $\Z-Vect$ are the same, except
we work in the category $\Z-Vect$. Specifically, a $\Z$-graded co-operad
is a collection of spaces $\mc(n)_k$, $k\in\Z$ and for $k=k_1+...+k_n$,
maps
\beg{eopgr}{\gamma:\mc(m_1+...+m_n)_k\r \mc(m_1)_{k_1}\otimes...\otimes \mc(m_n)_{k_n}}
which satisfy the obvious analogue of the axioms \rref{edop1}, \rref{edop2}.}
\end{definition}

\vspace{3mm}

\begin{remark}
\label{rwaste}
{\em The operations \rref{dopg1}
are symmetrical, but `wasteful' in the sense that they can be generated
by less general operations, called {\em insertions}, and
the notion of an operad can be equivalently formulated 
in terms of these operations. An insertion is
an operation 
\beg{erwaste}{\phi:\mc(m+1)\otimes \mc(n-m)\r \mc(n),
}
and is defined by
$$\phi(x,y)=\gamma(x,1,...,1,y).$$
The axioms are equivariance (the restriction of general equivariance 
to the cases which
apply to insertions), associativity of insertions, and a commutativity of
insertions into two different variables (insertion into a 
different variable than the last is defined using the
equivariance). This clearly extends to the other structures defined
in Definition \ref{dop}.
Both formulations are easily
seen to be equivalent. 
}
\end{remark}

\vspace{3mm}

\begin{remark}
\label{rconc}
{\em First note that there is a stronger variant of axiom
\rref{eopgr} in the case of graded co-operads, as \rref{eopgr} allows
an ``infinite sum'' when one varies the degrees $k_1,...,k_n$
with given sum $k$ on the right hand side. We could require that
all but finitely many terms in this sum are $0$. This stronger
axiom is {\em not} what we mean in this paper.

Next, note that $\Z$-graded operads  and co-operads model situations when we
have operations which shift total degree by a given number, but
the operations which apply to each tuple of degrees are exactly
the same. When we want to capture the situation where available operations
depend on the input degrees, we must introduce the notion
of {\em $\Z$-sorted operad}. This is a system of vector spaces
$V(m_0,m_1,...,m_n)$, $m_i\in\Z$ where the right action
of $\sigma\in \Sigma_n$ has
$$\sigma: V(m_0,m_1,...,m_n)\r V(m_0,m_{\sigma(1)},...,m_{\sigma(n)})
$$
and composition is defined
as 
$$\begin{array}{l}\gamma:V(m_0,m_1,...,m_n)\otimes
V(m_1,m_{11},...,m_{1,k_1})\otimes...\otimes V(m_n,m_{n1},...,m_{n,k_1})\r\\
V(m_0,m_{11},...,m_{nk_n}).\end{array}$$
The diagram axioms \rref{edop1}, \rref{edop2} are the same.
Because the axioims are again in shapes of diagrams,
we have a dual notion of $\Z$-sorted co-operads. In fact, here, one may
replace $\Z$ by any set $I$, and define $I$-sorted operads (or co-operads)
in the same fashion.}

\end{remark}

\vspace{3mm}

\begin{proposition}
\label{pap1}
The forgetful functor from the category of operads 
(and the obvious homomorphisms - maps preserving the operations) to vector spaces
is a right adjoint (similarly for the forgetful functor from graded and
$\Z$-sorted operads to the category of $\N\times\Z$-graded vector
spaces and vector spaces graded by $\coprod_{n\geq 0}\Z^n$).
The forgetful functor from the category of co-operads (resp. graded co-operads,
resp. $\Z$-sorted co-operads) to $\N$-graded vector
spaces (resp. the category of $\N\times\Z$-graded vector
resp. vector spaces graded by $\coprod_{n\geq 0}\Z^n$) is a left adjoint.
\end{proposition}

We shall postpone the proof because it is an example of
an even more general principle. 

\begin{definition}
\label{dalg}
{\em Let $\mc$ be an $I$-sorted operad. A {\em $\mc$-algebra}
is a system of vector spaces $X_i$, $i\in I$, and a system of
maps
\beg{ealg1}{\theta:\mc(i_0,i_1,...,i_n)\otimes X_{i_1}\otimes ...\otimes X_{i_n}
\r X_{i_0}
}
satisfying
\beg{ealgc}{(\theta\sigma)(?,...?)=\theta(\sigma(?,...,?))}
where $\theta\in\Sigma_n$ and $\sigma$ acts on tuples by permutation,
and
\beg{ealgc2}{\theta(\gamma(?,?,...,?),?,...,?)=\theta(?,\theta(?,?,...,?),...,\theta(?,?,...,?))
}
when applicable. }

\end{definition}

\begin{remark}
\label{ralg}
{\em In the previous definition, one may also define an endomorphism
operad and coendomorphism operad of $X$ by
\beg{eend}{End(X)(i_0,i_1,...,i_n)=Hom(X_{i_1}\otimes...\otimes X_{i_n}, X_{i_0}),
}
\beg{ecend}{Coend(X)(i_0,i_1,...,i_n)=Hom(X_{i_0},X_{i_1}\otimes...\otimes X_{i_n}).
}
Then a $\mc$-algebra is the same thing
as a homomorphism of $I$-sorted operads
\beg{eend1}{\mc\r End(X).}
It is therefore natural to define a coalgebra $X$ over an $I$-sorted operad $\mc$ to
be a morphism of operads
\beg{eend2}{\mc\r Coend(X).}
This, of course can also be written in terms of 
maps
\beg{ealg2}{\theta:\mc(i_0,i_1,...,i_n)\otimes X_{i_0}\r X_{i_1}\otimes ...\otimes X_{i_n}
}
and conditions analogous to \rref{ealgc}, \rref{ealgc2}. 

We can also define algebras resp. coalgebras over a co-operad
by dualizing \rref{ealg1}, \rref{ealgc}, \rref{ealgc2} resp. 
\rref{ealg2} and the corresponding conditions. However, in for these
dual notions,we do not know of an analogue of the descriptions \rref{eend1},
\rref{eend2}.}
\end{remark}

\vspace{3mm}

\begin{proposition}
\label{pop}
For every set $I$, there exists an $\coprod_{n\geq 0}I^{n+1}$-sorted operad
$Q_I$ such that the category of $I$-sorted operads (resp. co-operads) is canonically
equivalent to the category of $Q_I$-algebras (resp. $Q_I$-co-algebras).

\end{proposition}

\qed

\vspace{3mm}

Proposition \ref{pap1} now follows from the following result.

\vspace{3mm}

\begin{proposition}
\label{pap2}
Let $\mc$ be an $I$-sorted operad. Then the forgetful functor from
the category of $\mc$-algebras (resp. $\mc$-coalgebras) to $I$-graded vector spaces is
a right (resp. left) adjoint.
\end{proposition}

\Proof
The statement about $\mc$-algebras is classical. In effect,
the left adjoint is 
\beg{eopmm}{CX_i=\cform{\bigoplus}{n\geq 0}{}(\cform{\bigoplus}{(i,i_1,...,i_n)\in I^{n+1}}{} \mc(
i,i_1,...,i_n)\otimes X_{i_1}
\otimes...\otimes X_{i_n})/\Sigma_n,
}
and in fact the category of $\mc$-algebras is equivalent to the
category of algebras over the monad $C$. (For the definition of monads and
their algebras, see \cite{ml}.)

The statement about $\mc$-coalgebras follows from the fact that the forgetful
functor from $\mc$-coalgebras to $I$-graded vector spaces obviously
preserves coproducts and coequalizers, and one can use the adjoint
functor theorem to show it must have a right adjoint. 
\qed

\vspace{3mm}
\begin{remark}
\label{rcoalg}
{\em It is not true that the right adjoint of the forgetful functor from
$\mc$-coalgebras to $I$-graded vector spaces would be
given by a dual of \rref{eopm}. The dual is
\beg{eopm}{PX=\cform{\prod}{n\geq 0}{}(\cform{\prod}{(i, i_1,...,i_n)\in I^n}{} 
Hom(\mc(i,i_1,...,i_n),X_{i_1}
\otimes...\otimes X_{i_n})^{\Sigma_n},
}
but \rref{eopm} fails to be a comonad. The problem is that tensor product
does not distribute under the product, so when one writes the formula
for a map
$$PX\r PPX$$
dual to the classical monad structure on $C$, the image actually will be
a vector space containing but not equal to $PPX$. The actual comonad
is the intersection $P^\prime X$ of the inverse images of all such
maps $X\r P^n X$.}
\end{remark}

\vspace{3mm}

\begin{remark}
\label{rfail}
{\em Another example of failure of dualization is the fact that 
the forgetful functor from the category of algebras over an 
$I$-sorted cooperad
$\mc$ to the category of $I$-graded vector spaces does {\em not}
in general create products. Note that the structure map is
\beg{efail}{X_{i_1}\otimes...X_{i_n}\r \mc(i,i_1,...,i_n) \otimes X_i.
}
When we take a product of $\mc$-coalgebras
$X(j)$, 
then the image of the structure map \rref{efail} for $\prod X(j)$
ends up in 
$$\prod \mc(i,i_1,...,i_n)\otimes X(j)_i,$$
while $\mc$-coalgebra structure requires it to be in the subspace 
$$\mc(i,i_1,...,i_n)\otimes \prod X(j)_i.$$
However, note that when all of the spaces $\mc(i,i_1,...,i_n)$
are finite-dimensional, then the collection of the dual spaces $\mc(i,i_1,...,i_n)^{\vee}$
forms an $I$-sorted operad $\mc^{\vee}$, and the category of
$\mc$-algebras is equivalent to the categort of $\mc^{\vee}$-algebras,
which, as remarked above, is canonically equivalent to the
category of monads over the category of $I$-graded vector spaces.
In particular, in that case, the forgetful functor does 
create products.}
\end{remark}

\vspace{3mm}
Let $\mc$ be an $I$-sorted operad.
There is one more construction we need to cover, namely
quotient $\mc$-algebras, and the dual construction for
$\mc$-coalgebras. The most general version of the
story is this: Let 
\beg{eu}{U:C\r D
}
be a functor with left adjoint $L$, and let $X$ be an
object of $C$. Then we can consider the category $X/C$
of arrows $X\r ?$, where $?$ is an object or morphism
of $C$. The forgetful functor
$$X/C\r C$$
creates limits. Then $U$ induces a functor
\beg{eu1}{X/U:X/C\r UX/D,}
which therefore preserves limits, and by the adjunct functor
theorem has a left adjoint $X/L$. Similarly, when \rref{eu}
has a right adjoint $R$, and $X$ is an object of $C$, we may consider the
category $C/X$ of arrows $?\r X$, and the forgetful functor
\beg{eu2}{U/X:C/X\r D/UX}
preserves colimits, and hence has a right adjoint $R/X$ by
the adjoint functor theorem. Note that this construction
applies to the cases of $\mc$-algebras and $\mc$-coalgebras
by Proposition \ref{pap2}.

\vspace{3mm}

\begin{lemma}
\label{lfactor}
When $D$ is the category of $I$-graded vector spaces and
$C$ is the category of algebras (resp. coalgebras) over
an $I$-sorted operad, $U$ is the forgetful functor (see \rref{eu}) and
$L$ (resp. $R$) is the left (resp. right) adjoint, 
then the functor $X/L$ (resp. $R/X$) preserves epimorphisms
(resp. monomorphisms).
\end{lemma}

\Proof
The key point is that for a morphism
\beg{eab}{f:A\r B}
of $\mc$-algebras (resp. $\mc$-coalgebras),
the image of $f$ is a $\mc$-algebra (resp. $\mc$-coalgebra).
This factorizes any morphism \rref{eab} as 
\beg{egh}{\text{$f=gh$ where $g$ is mono and $h$ is epi.}} 
But if we take an epimorhpism (resp. monomorphism) $\phi$ in
$UX/D$ (resp. $D/UX$), and consider $f=X/L(\phi)$, then
factoring $f$ as in \rref{egh}, $h$ (resp. $g$) enjoys the
same universal property as $f$, and hence must be isomorphic to $f$.
\qed

\vspace{3mm}

\begin{remark}
\label{rfactor}
{\em In the case of $\mc$-algebras, we can interpret $X/L(\phi)$ where
$\phi$ is an epimorphism of $I$-graded vector spaces as
the quotient of $X$ by the ideal generated by $Ker(\phi)$. In
the case of $\mc$-coalgebras, the interpretation is
dual.

(Note that operads in our sense describe algebras with underlying
vector space structure and multilinear operations with multilinear
relations. In this setting, taking quotients with respect to
a congruence reduces to taking quotients by an ideal.)}
\end{remark}

\end{document}